\input amstex
\documentstyle{amsppt}
\vsize=8.0in \hsize 6.8 in
\voffset=2cm \loadbold \topmatter
\title On non-negatively curved metrics on open five-dimensional manifolds\endtitle
\author Valery Marenich and Mikael Bengtsson \endauthor
\thanks
\endthanks
\address H\"{o}gskolan i Kalmar, 391~82, Kalmar, Sweden\endaddress \email valery.marenich\@ hik.se, mikael.bengtsson\@ hik.se \endemail

\abstract
Let $V^n$ be an open manifold of non-negative sectional curvature with a soul $\Sigma$ of co-dimension two. The universal cover
$\tilde N$ of the unit normal bundle $N$ of the soul in such a manifold is isometric to the direct product $M^{n-2}\times R$.
In the study of the metric structure of $V^n$ an important role plays the vector field $X$ which belongs to the projection of
the vertical planes distribution of the Riemannian submersion $\pi:V\to\Sigma$ on the factor $M$ in this metric splitting
$\tilde N=M\times R$. The case $n=4$ was considered in [GT] where the authors prove that $X$ is a Killing vector field while
the manifold $V^4$ is isometric to the quotient of $M^2\times (R^2,g_F)\times R$ by the flow along the corresponding Killing
field. Following an approach of [GT] we consider the next case $n=5$ and obtain the same result under the assumption that the
set of zeros of $X$ is not empty. Under this assumption we prove that both $M^3$ and $\Sigma^3$ admit an open-book
decomposition with a bending which is a closed geodesic and pages which are totally geodesic two-spheres, the vector field $X$
is Killing, while the whole manifold $V^5$ is isometric to the quotient of $M^3\times (R^2,g_F)\times R$ by the flow along
corresponding Killing field.
\endabstract
\subjclass   53C20, 53C21. Supported by the Faculty of Natural Sciences of the Hogskolan i Kalmar, (Sweden)
\endsubjclass
 \dedicatory  \enddedicatory
\endtopmatter

\document
\baselineskip 12pt

\head \bf 1. Introduction \endhead

Let $(V^{n},g)$ be a complete open Riemannian manifold of non-negative sectional curvature. Remind that as follows from [CG]
and [P] an arbitrary complete open manifold $V^n$ of non-negative sectional curvature contains a closed absolutely convex and
totally geodesic submanifold $\Sigma$ (called a soul) such that the projection $\pi: V\to \Sigma$ of $V$ onto $\Sigma$ along
geodesics normal to $\Sigma$ is well-defined and is a Riemannian submersion (see also [CaS]). The (vertical) fibers
$F_P=\pi^{-1}(P), P\in \Sigma$ of $\pi$ define a metric foliation in $V$ and two distributions: a vertical ${\Cal V}$
distribution of subspaces tangent to fibers and a horizontal distribution ${\Cal H}$ of subspaces normal to ${\Cal V}$. For an
arbitrary point $P$ on $\Sigma$, an arbitrary geodesic $\gamma(t)$ on $\Sigma$ and arbitrary vector field $V(t)$ which is
parallel along $\gamma$ and normal to $\Sigma$ the following
$$
\Pi(t,s)=exp_{\gamma(t)} sV(t) \tag 1
$$
are totally geodesic surfaces in $V^n$ of zero curvature, i.e., flats.

\medskip

When $dim \Sigma=1$ or $codim(\Sigma)=1$ the manifold $V^n$ is locally isometric to the direct product of $\Sigma$ and
Euclidean space of a complementary dimension and of non-negative curvature. Study of the next case $codim(\Sigma)=2$ was began
in [M1], where we noted that the manifold $V^4$ or is a direct product when the holonomy of the normal bundle of $\Sigma$ in
$V$ is trivial, or the holonomy group acts transitively on normal vectors, every geodesic normal to $\Sigma$ is a ray and (1)
holds. The metric structure in this case might be more complicated. In [GT] the authors consider four-dimensional manifolds
diffeomorphic to direct products $M^2\times R^2$ and prove the following.

\medskip

\proclaim{Theorem~A, [GT]}Every non-negatively curved metric on $M^2\times R^2$ is isometric to a Riemannian quotient of the
form $((M^2,g_0)\times(R^2,g_F)\times R)/R$. Here $R$ acts diagonally on the product by the flow along Killing vector fields on
$(M^2,g_0)$ and $(R^2,g_F)$ and by translations on $R$.
\endproclaim

\medskip

The very important role in the proof of the Theorem~A plays the vector field $X$ which is the projection of the vertical vector
field in the universal cover $\tilde N$ of the boundary of some metric $s$-tube $N$ of the soul on the "horizontal" factor $M$
in the metric splitting $\tilde N=M\times R$, see below. In the four-dimensional case this vector field $X$ restricted to $M^2$
always has zeros since $M^2$ is a two-dimensional sphere. In our case the soul $\Sigma$ of $V$ and $M^3$ are three-dimensional
spheres, and hence, $X$ might be nowhere zero as the following simple example shows. Let $h:S^3\to S^2$ be the Hopf bundle,
i.e., the factoring of a unit sphere $S^3$ in the complex plane $C^2$ by $S^1$ action - multiplication by complex numbers of
absolute value $1$. Consider $V^5$ which is the the quotient of the direct product of $S^3\times R^2\times S^1$, where the
$S^1$ acts on $R^2=C$ by rotations, i.e., again, multiplication in $C^1$ by unit complex numbers. Then for the manifold
$V^5=S^3\times R^2\times S^1 / S^1$ the vector field $X$ is nowhere zero.\footnote{For corresponding $M^3$ and $S^3$ the
one-form given by the scalar product with $X$ is a (nowhere degenerated) contact form $\alpha$ with $\alpha\wedge d\alpha$ -
the volume form.} The objective of this note is to expand an approach from [GT] to the case of non-negatively curved
five-dimensional $V^5$ diffeomorphic to a direct product $S^3\times R^2$ under the following assumption.

\medskip

{\bf Assumption~1. The set of zeros of the vector field $X$ is not empty.}

\medskip

Our main result is very similar to the Theorem~A above.

\proclaim{Theorem~B} Let $V^5$ be an open manifold of non-negative sectional curvature and difeomorphic to $S^3\times R^2$.
Assume that the vector field $X$ has non-empty zero set ${\Cal Z}$. Then ${\Cal Z}$ is a closed geodesic and the manifold $M^3$
admits a singular foliation - "open-book decompositions" by totally geodesic and isometric to each other horizontal
two-dimensional spheres $S^2(\psi)$, where the singular set of this decompositions - "bindings", equal the closed geodesic
${\Cal Z}$. The flow along Killing field $X$ acts as "turning pages" in this open-book decomposition, while $V^5$ itself is
isometric to a Riemannian quotient of the form $M^3\times(R^2,g_F)\times R/R$ with $R$ acting diagonally on the product by the
flow along Killing vector fields on $M^3$ and $(R^2,g_F)$ and by translations on $R$. The Riemannian submersion
$\pi:V^5\to\Sigma$ conveys the open-book decomposition of $M^3$ to a similar open-book decomposition of $\Sigma$ with the pages
$\Sigma^2(\psi)$ isometric to $S^2(\psi)$.
\endproclaim

\medskip

In the same way as Theorem~A in [GT] our Theorem~B follows from the fact that the vector field $X$ on $N$ is Killing for every
$s$, where $N$ is the boundary of $s$-metric neighborhood of the soul $\Sigma$ in $V$, see Theorem~3 below. Thus, after proving
Theorem~3, we complete the proof of Theorem~B by referring to the corresponding arguments from [GT], see section~5.

Note that the general case of five-dimensional open manifold $V^5$ with a soul of codimension $2$ can be reduced to the one
under consideration as follows. First, we note that if the fundamental group of $\Sigma$ (which is isomorphic to that of $V$)
is not finite, the universal cover $\tilde V$ contains a straight line in the universal cover $\tilde\Sigma$ of the soul. Then
both $\tilde V$ and $\tilde\Sigma$ split into direct products, and the case is reduced to the already studied one of open
four-dimensional manifolds. When the fundamental group of $\Sigma$ is finite the universal cover $\tilde\Sigma$ is
diffeomorphic to a sphere $S^3$ due to the non-negativity of the curvature. Next: because an arbitrary vector bundle over
simply connected $S^3$ is, obviously, trivial we see that an investigation of the metric structure of an arbitrary $V^5$ with a
soul of codimension $2$ is reduced to the case when $V^5$ is diffeomorphic to the direct product $S^3\times R^2$.\footnote{The
case when five-dimensional $V^5$ has a soul of codimension $3$ we considered in [M5].}

Below we assume that the holonomy of the normal bundle is not trivial, for otherwise by a direct product theorem from [M1,4]
the manifold $V$ is a metric product.

\head \bf 2. Vector field $X$ and its zeros\endhead

Fix some positive $s_0$ smaller than a focal radius of $\Sigma$ in $V$. For some $s<s_0$ denote by $N\Sigma(s)$, or simply by
$N$, the boundary of an $s$-neighborhood of $\Sigma$. Due to our choice it is a smooth manifold. It consists of all points
$Q(P,V)=exp_P(sV)$, where $P$ is a point on $\Sigma$ and $V$ is a unit vector normal to $\Sigma$ at $P$.

\medskip

\proclaim{Lemma~1} $N(s)$ has non-negative curvature if $s$ is sufficiently small. \endproclaim

\medskip

\demo{Proof} This follows from the Gauss equations and the fact that $N(s)$ bounds a convex subset in a manifold $V$ of
non-negative curvature. The last is obviously true when the holonomy of the (trivial) normal bundle $\nu\Sigma$ of the soul is
trivial, i.e., all parallel translations along closed curves in $\Sigma$ acts identically  on vectors normal to $\Sigma$
because then $V$ is isometric to the direct metric product $\Sigma\times(R^2,h)$, see [M1]. If the holonomy is not trivial,
then all normal vectors are so called ray directions, and $N(s)$ coincides with the boundary $\partial C_s$ of an absolutely
convex set constructed in [CG], see again [M1]. The Lemma~1 is proved.
\enddemo

\medskip

\proclaim{Lemma~2} The universal cover $\tilde N$ of $N(s)$ is isometric to the direct product $(M,g)\times R$, where $M$ is
diffeomorphic to $S^3$. The composition of a covering map and a submersion $\pi$ provides a diffeomorphism between an arbitrary
factor $M$ and the soul $\Sigma$ which we denote by $\pi^M: M\to \Sigma$.\footnote{Note, that this statement and forthcoming
(3) both are true for an arbitrary $V^n$ with simply connected soul of codimension two.}
\endproclaim

\medskip

\demo{Proof} This follows from the fact that $N(s)$ is diffeomorphic to the trivial circle bundle over three-dimensional sphere
$\Sigma$, i.e., has an infinite cycle fundamental group generated by a homotopy class of a fiber. Then by standard arguments
the universal cover $\tilde N(s)$ admits a straight line, and hence by Toponogov splitting theorem is isometric to the direct
product $(M,g)\times R$.

Denote by $E$ the unit vector field in $N(s)$ tangent to the projections of straight lines (i.e., $R$-factor) from $\tilde N$
to $N$. By $W$ we denote the (vertical) vector field on $N$ which is the speed of the natural $S^1$-action on $N$ given by
rotations in a positive direction of a normal vectors to $\Sigma$ as follows: for $Q=Q(P,V)$ denote by
$Q_\phi=Q_\phi(P,V)=Q(P,V_\phi)$, where $V_\phi$ is $V$ rotated by the angle $\phi$ in the bundle of unit normals to $\Sigma$
in $V$ (which is correctly defined since the bundle is topologically trivial). Finely, denote by $X$ the vector field on $N$
which is the component of $W$ normal to $E$.\footnote{Note, that our $X$ is different from similar $X$ of [GT].}
$$
X=W-(W,E)E. \tag 2
$$
Note, that $N$ naturally inherits from $V$ a horizontal distribution ${\Cal H}$, while the vector field $W$ belongs to the
vertical distribution. If by $M$ we denote an image of some $(M,g)$-factor in the direct metric product $\tilde N$ under the
projection $pr:\tilde N\to N$, then (the restriction of) $E$ on $M$ would be the unit vector field of normals to $M$, $X$ is a
vector field tangent to $M$, while another vector field $Y$ tangent to $M$ would be a horizontal if and only if it is normal to
$X$. In particular, the tangent subspace $T_QM$ is horizontal ${\Cal H}_Q$ if and only if $X(Q)=0$.

Note that the vector field $W$ in $N$ is never tangent to any of the $M$-factor, or (equivalently) never orthogonal to $E$.
Indeed, if so then some homotopicaly non-trivial closed geodesic $\Gamma(s)$ in $N$ which is the images of a straight line in
the universal cover $\tilde N=M\times R$, would be horizontal at some point, and therefore, horizontal everywhere, which
obviously can not be homotopicaly non-trivial in $N$. To see this denote by $\bar\Gamma(s)$ its image under $\pi$ in $\Sigma$,
and by $V(s)$ the normal vector field of vertical geodesics connecting $\bar\Gamma(s)$ and $\Gamma(s)$. Since $\Sigma$ is
simply connected their exists a disk $D$ in $\Sigma$ with a boundary $\Gamma$ and extension of the vector field $V$ over $D$
(because the restriction of a normal bundle to $D$ is trivial). The vertical lift of $D$ along this extension will provide us a
disk in $N$ with a boundary $\Gamma$ implying that $\Gamma$ is contractible in $N$. The obtained contradiction proves that $E$
is never horizontal, or that the map $\pi^M:M\to\Sigma$ from any of the image $M$ of a factor in the direct product $\tilde
N=(M^3,g)\times R$ into the soul $\Sigma$ is a diffeomorphism. The Lemma~2 is proved.\footnote{This also fills the gap in the
arguments from Lemma~2.1 in [GT].}
\enddemo

\medskip

By definition the differential of the diffeomorphism $\pi^M: M\to \Sigma$ is an isometry on the subspace of horizontal vectors,
i.e., on the subspace in $T_QM$ normal to $X$ (or on the whole $T_QM$ if $X(Q)=0$), while
$$
\|d\pi^M_Q(X)\|=cos(\alpha(Q))\|X\|, \tag 3
$$
where $\alpha(Q)$ denotes the angle at the point $Q\in M$ between vectors $E$ and $W$. The map $\pi^N:N\to\Sigma$ is the
composition of the projection in the universal cover to the horizontal factor and then $\pi^M$.

When $X$ is identically zero the submanifold $M$ is horizontal in $N$, isometric to $\Sigma$ by (3), the holonomy of the normal
bundle $\nu\Sigma$ is trivial, and, again, $V$ is isometric to the direct product $\Sigma\times (R^2,h)$ of the soul $\Sigma$
and some non-negatively curved plane $(R^2,h)$.

Next we prove that if $X$ is not identically zero, or has no zeros at all, then $X$ vanish along some closed geodesic.

\medskip

\proclaim{Theorem~1} If the set of zeros ${\Cal Z}$ of the vector field $X$ in $M$ is a proper subset (i.e., is not $M$ itself
or empty) then ${\Cal Z}$ is a closed geodesic. Every minimal geodesic connecting two points from ${\Cal Z}$ is itself a subset
of ${\Cal Z}$.
\endproclaim

\medskip

\demo{Proof} If ${\Cal Z}$ is a proper subset of $M$ then for some $P\in{\Cal Z}$ there exists a sequence of points $Q_i\to P$
such that $X(Q_i)\not= 0$. As in [GT], see Lemma~2.1;  we note that every geodesic $L(P,Q;t)$ in $M$ connecting a point $P$
where $X$ vanish with an arbitrary point $Q$ with non-vanishing $X(Q)$ is orthogonal to $X(Q)$,
$$
\bar{QP} \perp X(Q), \tag 4
$$
where $\bar{PQ}$ denotes the vector of direction of $L(P,Q;t)$ at the point $Q$. Hence, ${\Cal Z}$ belongs to the exponential
image $\Pi(Q,X(Q))$ of a plane in $T_QM$ of all vectors normal to $X(Q)$:
$$
{\Cal Z} \subset \Pi(Q,X(Q)). \tag 5
$$
The surface $\Pi(Q,X(Q))$ near $Q$ is "almost a plane" - smooth with a second form vanishing at $Q$. Fix for a moment some
$Q=Q_i$ close enough to $P$. Then in a small closed ball $B$ around $P$ with radius $dist(P,Q)$ zeros of $X$ belong to this
"almost a plane" $\Pi(Q,X(Q))$. Thus there exists the farthest point $Q'\in B$ to $\Pi$ where $X(Q')\not= 0$. Then by (5) the
part of the set ${\Cal Z}$ inside the ball $B$ belongs to the intersection of two "almost planes" $\Pi(Q,X(Q))$ and
$\Pi(Q',X(Q'))$ which are "almost orthogonal". This intersection, as easy to see, is a smooth curve with geodesic curvature of
the order $dist(P,Q)$. Because the point $Q=Q_i$ can be chosen arbitrary close to $P$ we conclude that the set ${\Cal Z}$ is
inside some finite collection of intervals of geodesics. Next, we verify that ${\Cal Z}$ is connected. Indeed, if not we may
find two different points $P_1$ and $P_2$ from its different components and such that the minimal geodesic $L=L(P_1,P_2; t)$
connecting these points does not intersect ${\Cal Z}$. In some small ball $B$ around the middle point $Q$ of this geodesic the
vector field $X$ will be non-zero with $X(Q)$ normal to $L$. Consider the "almost plane" $\Pi$ in $B$ going through the point
$Q$ and normal to $L$. From (4) we see that the vector field $X$ in this plane not only is almost tangent to $\Pi$, but also
almost tangent to small circles in $\Pi$ around $Q$. Which implies that the projection of $X$ on $\Pi$ has index $\pm 1$ at the
center $Q$ of these circles, i.e., equals zero at $Q$. The obtained contradiction proves that ${\Cal Z}$ is a closed
geodesic.\footnote{By the arguments above we immediately deduce that ${\Cal Z}$ is a connected geodesic. The fact that this
geodesic can not be infinite follows from the compactness of $M$ and that it can not accumulate to something other than
itself.} Clearly, if some minimal geodesic $L$ connects two zeros $P$ and $Q$ from ${\Cal Z}$, but does not belong to ${\Cal
Z}$ we may repeat arguments above to show that there exists one more point in the interior of $L$ where $X$ vanish. Which
completes the proof of our theorem.
\enddemo

\medskip

Note, that in our arguments we used only condition (4). Hence we have the following.

\medskip

\proclaim{Corollary~1} If some smooth (not identically zero) vector field $X$ in some compact three-dimensional manifold $M$
satisfies (4), then its set of zeros ${\Cal Z}$ is a closed geodesic.
\endproclaim

\medskip

Next we consider the metric structures of $M$ when ${\Cal Z}$ is a closed geodesic.\footnote{When $X$ has no zeros we may
introduce the following {\bf $A$-"contact" structure} on $M$. Denote by $\alpha$ the 1-form on $M^3$ given by the scalar
product with a vector filed $X$. Because the vector field $E$ on $N$ is parallel, from $\nabla^N E\equiv 0$ (here $\nabla^N$ is
the covariant derivative in $N$ in a metric induced by $N\subset V$) we see that
$$
(\nabla^N_Y X, Z)=(\nabla_Y W - (\nabla_Y W),E)E,Z)=(\nabla_Y W,Z)
$$
for arbitrary $Y,Z$ tangent to $M$. Therefore, because $M$ is totally geodesic in $N$ we have by direct calculations that
$$
\alpha\wedge d\alpha = a(P) d vol^M,
$$
where $d vol^M$ is the volume form of $M$ and the function $a(P)$ is given by
$$
a(P)=(A_Y Z, X)
$$
or by $(A_Y' Z', X')$ where $\{X',Y',Z'\}$ an arbitrary orthonormal (positively orientated) basis in $T_P M$. The horizontal
distribution on $M$ is not involutive outside zeros of $a$. It would be interesting to find examples with $a$ vanishing
somewhere and nowhere zero $\alpha$.}

\head \bf 3. ${\Cal Z}=S^1$ \endhead

Assume that the set of zeros ${\Cal Z}$ of the vector field $X$ is some closed geodesic ${\Cal Z}={\Cal Z}(t), 0\leq t\leq 1$.
Then ${\Cal Z}(t)$ is horizontal, its projection by the submersion $\pi$ to the soul $\Sigma$ is again a closed geodesic $\bar
{\Cal Z}(t)$ of the same length, for every $t$ the geodesic $l_t(s)$ connecting $\bar {\Cal Z}(t)$ and ${\Cal Z}(t)$ is normal
to $\Sigma$ with a direction $V(t)$ parallel along $\bar {\Cal Z}(t)$. Now take an arbitrary point $Q$ in $M$ and connect it
with all the points ${\Cal Z}(t)$ by minimal geodesics $L_t(s)$. All this geodesics are horizontal, and, if $\bar Q=\pi(Q)$,
their projections $\bar L_t(s)$ are minimal geodesics connecting $\bar Q$ with $\bar {\Cal Z}(t)$. Also, if $V(t,s)$ denotes
the (unit) vector of the direction of the (vertical) geodesic connecting $\bar L_t(s)$ with $L_t(s)$, then $V(t,s)$ is parallel
along $\bar L_t(s)$. In particular, it follows that the parallel translation along a closed path from $\bar Q$ to $\bar Q$
consisting of two $\bar L_{t'}(s)$ and $\bar L_{t"}(s)$ and a part $\bar L(t), t'\leq t\leq t"$ acts trivially on $V$ - the
direction of $\bar Q Q$. Which by the prism construction from [M1-3] implies that the O'Neill's fundamental tensor $A$ vanishes
at $Q$ for horizontal vectors tangent to the family of geodesics $L_t(s)$. As we already saw, this family belongs to the plane
$\Pi(Q,X(Q))$ of all geodesics, issuing from $Q$ in directions normal to $X(Q)$. Thus we have
$$
A_Y Z(Q)\equiv 0 \tag 6
$$
for all $Y,Z$. Again, by the same prism construction we have
$$
R(\bar Y(t),\bar Z(t))V(t) \equiv 0 \tag 7
$$
along $\bar {\Cal Z}(t)$, where $\bar Z(t)$ is the unit tangent to $\bar {\Cal Z}(t)$, $\bar Y(t)$ any tangent to $\Sigma$, and
$V(t)$ is the direction of the vertical geodesic $\bar {\Cal Z}(t){\Cal Z}(t)$.

Because $Q$ (and $\bar Q$ correspondingly) was arbitrary, the tensor $A$ vanishes identically in $M$ on vectors normal to the
vector field $X$.

\medskip

\proclaim{Theorem~2} Distribution in $M\backslash{\Cal Z}$ of the two-planes normal to the vector field $X$ is integrable. It
is tangent to the family of totally geodesic spheres with a common intersection set - the closed geodesic ${\Cal Z}$.
\endproclaim

\medskip

\demo{Proof} Indeed, from (6) immediately follows that the Lie bracket of arbitrary fields $Y,Z$ orthogonal to $X$ is also
orthogonal to $X$:
$$
([Y,Z], X)=([Y,Z], W-(W,E)E)=([Y,Z],W)=(A_YZ-A_ZY)=0. \tag 8
$$
Thus, the vector field $X$ is, actually, the field of normals to some family of hyper-surfaces in $M$. But, as we already know,
every geodesic $L_t$ connecting $Q$ with a point ${\Cal Z}(t)$ of ${\Cal Z}$ belongs to such a surface. From which, obviously,
follows that this family of surfaces coincide with the family of our planes $\Pi(Q,X(Q))$. Because at $Q$ the second form of
this surface vanish, and $Q$ is arbitrary our surfaces have vanishing second forms or are totally geodesic. For definiteness,
from now on we call by $\Pi(Q)$ the union of all geodesics $L_t$ connecting $Q$ with ${\Cal Z}$. It is a totally geodesic
surface which boundary is the closed geodesic ${\Cal Z}(t)$. Therefore, the vector $Y(t)$ tangent to $\Pi$ and normal to this
boundary at the point ${\Cal Z}(t)$ is parallel along ${\Cal Z}$. Because the tangent vector $Z(t)$ to this geodesic is also
(auto-)parallel, we see that the holonomy around ${\Cal Z}$ is trivial, i.e., parallel translation along ${\Cal Z}(t)$ is the
identity operator. If we choose some parallel vector field $Y^*(y)$ along ${\Cal Z}$ normal to $Z(t)$, we can define the angle
function $\psi$ for vectors $Y(t)$ normal to $Z(t)$ as the the angle between $Y(t)$ and $Y^*(t)$. Corresponding $\Pi(Q)$ we
denote also by $\Pi(\psi)$. To complete the proof of the theorem we note that for (a half-sphere) $\Pi(\psi)$ there exists
another one $\Pi(\psi+\pi)$ which normal to their common boundary ${\Cal Z}(t)$ equals $-Y(t)$. Their union in a neighborhood
of ${\Cal Z}(t)$ is again an exponential image of planes tangent to $Y(t)$ and $Z(t)$, and therefore, is a smooth surface: a
sphere which we denote by $S^2(Q)$, or by $S^2(\psi)$ (then $S^2(\psi)=S^2(\psi+\pi)$). Theorem~2 is proved.
\enddemo

\medskip

Configuration we described in the last theorem is well-known and is called {\bf an open book decomposition}.

\medskip

\proclaim{Corollary~2} If the set of zeros of $X$ is a closed geodesic ${\Cal Z(t)}$, then $M$ admits an open book
decomposition with a bending ${\Cal Z}$ and pages $\Pi(\psi)$ which are totally geodesic half-spheres.
\endproclaim

\medskip

Next we look more closely on the family of diffeomorphisms between pages $\Pi(\psi)$ of our open book decomposition given by
shifts in directions normal to them. Let $f_\theta :\Pi(\psi)\to \Pi(\psi+\theta)$ denotes the map sending the point $Q$ in
$\Pi(\psi)$ into the intersection of $\Pi(\psi+\theta)$ with an integral curve of the field of normals to pages issuing from
$Q$. If we denote
$$
\partial f_\theta (Q)/\partial\theta=X^*(f_\theta(Q)), \tag 9
$$
then the field $X^*$ is proportional to $X$, i.e., $X^*(Q)=k'(Q)X(Q)$ for some positive function on $M\backslash{\Cal Z}$. By
$k$ we denote its norm: $k(Q)=\|X^*(Q)\|$. Because all pages $\Pi(\psi)$ are totaly geodesic all maps $f_\theta$ are
isometries. Therefore, we call the family of these isometries: "turning pages". If, in addition, $k$ is constant along
trajectories of $X^*$ (or $X$, which is the same) then $f_\theta: M\to M$ is a family of isometries of the entire
$M$,\footnote{"reading the book"} and the vector field $X^*$ is a Killing vector field. Note also the following trivial
statement.

\medskip

\proclaim{Lemma~3} All trajectories of the vector field $X^*$ in $M$ are closed circles around ${\Cal Z}$.
\endproclaim

\medskip

\demo{Proof} Indeed, take some geodesic $L_t(s)$ in $S^2(\psi)$ connecting some $Q$ with the point ${\Cal Z}(t)$ which is
nearest to it, and consider the orbit of this geodesic under our family of "rotations": $\Phi(s,\theta)=f_\theta(L_t(s))$. We
choose natural parameter $s$ on $L_t(s)$ in such a way that ${\Cal Z}(t)=L_t(0)$. Because for every $\theta$ the curve
$f_\theta(L_t(s))$ lies in the totally geodesic $\Pi(\psi+\theta)$ and $f_\theta$ is an isometry, this $f_\theta(L_t(s))$ is
again the geodesic in $M$. Therefore, $\Phi(s,\theta)$ is a part of the "plane" $\Pi({\Cal Z}(t), Z(t))$ of all geodesics
issuing from ${\Cal Z}(t)$ in directions normal to $Z(t)$. For a fixed $s$ the line $\Phi(s,\theta)$ is a closed circle in this
"plane".\footnote{It is interesting to note also the following property of these circles. As we will show, the vector field
$X^*$ is Killing and constant along its trajectories, i.e., circles $f_\theta(Q)$. Therefore, the norm $k$ of $X^*$ attains its
maximum on some set ${\Cal Z}^*$ which is invariant under rotations $f_\theta$. We claim that ${\Cal Z}^*$ is a collection of
closed geodesics in $M$. Indeed, from $Y(X^*,X^*)\equiv 0$ for every $Y$ tangent to $S^2(\psi)$ at some point $Q$ of ${\Cal
Z}^*$ it follows that $\nabla_{X^*}X^*\equiv 0$, or that the geodesic curvature of the orbit $f_\theta(Q)$ equals zero. Every
closed geodesic from ${\Cal Z}^*$ is linked with ${\Cal Z}$.}
\enddemo

\medskip

\head \bf 4. $X^*$ is Killing \endhead

As we saw above, $A_Z$ vanishes along ${\Cal Z}$, see (6,7). Also the holonomy of the normal bundle is trivial along the
projection  $\bar {\Cal Z}(t)$ of ${\Cal Z}$ under submersion $\pi$ which is the closed geodesic in the soul $\Sigma$.
Therefore, applying the simplified version of arguments\footnote{when $codim \Sigma=2$ instead of $codim \Sigma=3$.} from the
proof of the Theorem~A from [M5] (see section~5 there) we get
$$
\nabla_W W\equiv 0 \qquad \text{ and } \qquad R[W(t),Z(t)]\equiv 0 \tag 10
$$
for a unit vertical field $W(t)$ along ${\Cal Z}$.

Take another vector $V(\phi)$ normal to $\Sigma$ at $\bar P=\bar {\Cal Z}(0)$ with an angle $\phi$ to $V$. Its parallel
transport along $\bar {\Cal Z}(t)$ is again $V(\phi)$. Denote by $V(\phi,t)$ the corresponding parallel vector field along
$\bar {\Cal Z}(t)$. The vertical lifts of $\bar {\Cal Z}(t)$ into $N$ along this vector field are again closed geodesics which
we denote by ${\Cal Z}(\phi,t)$. Easy to see that (6,7) are satisfied along them\footnote{for the proof note, that (7) implies
(6) through the prism construction, see the Lemma~2 in [M5]}, which in turn implies (10) along ${\Cal Z}(\phi,t)$, or that the
vertical fibers of the submersion $\pi: {\Cal Z}(\phi,t)\to \bar {\Cal Z}(t)$ have zero geodesic curvature, or are geodesic
lines in $N$. Hence, they coincide with projections of straight lines, i.e., $R$-factors under universal cover $M\times R\to
N$. We formulate the obtained result as follows.

\medskip

\proclaim{Lemma~4} The set of zeros of the vector field $X$ in $N$ is a tori which is the image of the direct product of ${\Cal
Z}\subset M$ with a straight-line factor $R$ in the universal cover $\tilde N=M\times R$ under covering map $\tilde N\to N$.
For an arbitrary choice of $M$ in $N$ the $\pi$-projection of the set of zeros of $X$ in $M$ is the same closed geodesic
$\bar{\Cal Z}$ in $\Sigma$.
\endproclaim

\medskip

The obtained claim means that every vertical fibre in $N$ stays in the set of zeros of the vector field $X$ if it contains some
of the point where $X$ vanish. Now we can repeat arguments from [GT] and prove the following statement.

\medskip

\proclaim{Theorem~3} The vector field $X^*$ is Killing, if it has non-empty zero set.
\endproclaim

\medskip

\demo{Proof} Indeed, the Lemma~4's claim enable us to repeat arguments from [GT]: for every point $Q$ of $M$ denote by
$(f_t(Q),t)$ the points of the fiber of the submersion $\pi:N\to\Sigma$ issuing from $Q$. These are trajectories of the vector
field $W$ in $N$. As we saw, the distance between $f_t(Q)$ and $f_t(P)$ is constant for every $P$ from the zero set ${\Cal Z}$
since the geodesic connecting them in $M\times\{t\}$ is horizontal. By the Lemma~4 we see $f_t$ is the identity map on ${\Cal
Z}$. Therefore, $f_t(Q)$ are circles $S^1_Q$ around ${\Cal Z}$, they coincide with the circles $f_\theta(Q)$ which are orbits
of the vector field $X^*$ above. To show that $X^*$ is Killing consider the cylinder $C_Q=\{(f_\theta(Q),t)\}$, (see Lemma~2.1
in [GT]). The restriction of $\pi$ on $C_Q$ is a Riemannian submersion of a flat cylinder onto some circle in $\Sigma$, or by
[GG] has fibers tangent to some Killing field on $C_Q$. This proves that $X^*$ has constant norm along $C_Q$ and is a Killing
vector field.
\enddemo

\medskip

\head \bf 5. Proof of the Theorem~B\endhead

From Theorem~3 it follows that the restriction of the Riemannian submersion $\pi:V\to\Sigma$ on $N$, which is the boundary of
some $s$-metric neighborhood of the soul, can be described as the factoring by the action along trajectories of the Killing
vector field $X^*$. From this fact the Theorem~B follows in the same way as Theorem~A; see section~3 in [GT] for the meticulous
analysis of the cooperation between Killing vector fields $X^*$ on different $s$-metric neighborhoods of the soul which ensures
the claim of both Theorems~A and B.

\enddocument
\bye